\documentclass{1}
\usepackage{amssymb,amsmath}
\usepackage[english,francais]{babel}
\newtheorem{theorem}{Theorem}%[section]

\newtheorem{e-proposition}[theorem]{Proposition}

\newtheorem{e-definition}[theorem]{Definition\rm}

\renewcommand{\theequation}{\arabic{equation}}
\def\og{\leavevmode\raise.3ex\hbox{$\scriptscriptstyle\langle\!\langle$~}}
\def\fg{\leavevmode\raise.3ex\hbox{~$\!\scriptscriptstyle\,\rangle\!\rangle$}}

\newcommand{\C}{{\mathbb{C}}}
\newcommand{\G}{{\mathbb{G}}}
\renewcommand{\P}{{\mathbb{P}}}
\newcommand{\R}{{\mathbb{R}}}

\newcommand{\calO}{{\mathcal{O}}}

\renewcommand{\div}{{\rm div}}
\DeclareMathOperator{\Br}{Br}
\DeclareMathOperator{\Gal}{Gal}
\DeclareMathOperator{\Pic}{Pic}

\DeclareMathOperator{\Spec}{Spec}

\newcommand{\wt}[1]{\widetilde{#1}}
\newcommand{\labelto}[1]{\overset{#1}{\longrightarrow}}

\newcommand{\ol}[1]{\overline{#1}}
\newcommand{\sing}{{\rm sing}}

\DeclareTextFontCommand{\textnf}{\normalfont}
\newcommand{\et}{{\textnf{\'et}}}

\begin{document}
\mbox{\ }\vspace{-105pt}

\begin{frontmatter}

\selectlanguage{english}
\title{A New Proof of Hilbert's Theorem on Ternary Quartics}
%\selectlanguage{francais}
%\title{ Une nouvelle d\'emonstration du th\'eoreme de Hilbert sur les quartiques
%ternaires.}

%%%%%%%%%%%%%%%%%%%%%%%%%%%%%%%%%%%%%%%%%%%%%%%%%%%%%%%%%%%%%%%%%%%
\selectlanguage{english}
\author[VP]{Victoria Powers}%\thanksref{labelVP}}
\ead{vicki@mathcs.emory.edu}
\address[VP]{Department of Mathematics and Computer Science,
          Emory University, Atlanta, GA 30322, USA}

%%%%%%%%%%%%%%%%%%%%%%%%%%%%%%%%%%%%%%%%%%%%%%%%%%%%%%%%%%%%%%%%%%%
\author[BR]{Bruce Reznick}%\thanksref{labelBR}}
\ead{reznick@math.uiuc.edu}
\address[BR]{Department of Mathematics, University of Illinois, Urbana, IL 61801, USA}

%%%%%%%%%%%%%%%%%%%%%%%%%%%%%%%%%%%%%%%%%%%%%%%%%%%%%%%%%%%%%%%%%%%
\author[CS]{Claus Scheiderer\thanksref{labelCS}}
\ead{claus@math.uni-duisburg.de}
\address[CS]{ Fachbereich f\"ur Mathematik und Statistik, Universit\"at 
Konstanz, 78457 Konstanz, Allemagne  }

\thanks[labelCS]{supported by European RTN-Network
   HPRN-CT-2001-00271 (RAAG).}

%%%%%%%%%%%%%%%%%%%%%%%%%%%%%%%%%%%%%%%%%%%%%%%%%%%%%%%%%%%%%%%%%%%
\author[FS]{Frank Sottile\thanksref{labelFS}}
\address[FS]{Department of Mathematics, Texas A\&M University,
             College Station,  Texas 77843,  USA} 
\ead{sottile@math.tamu.edu}
 \thanks[labelFS]{supported by the Clay Mathematical Institute,
  NSF CAREER grant DMS-0134860, and the MSRI.}
%\urladdr{http://www.math.tamu.edu/\~{}sottile}

%%%%%%%%%%%%%%%%%%%%%%%%%%%%%%%%%%%%%%%%%%%%%%%%%%%%%%%%%%%%%%%%%%%%%
%\mbox{\ }\vspace{-25pt}

\selectlanguage{english}
\begin{abstract}
Hilbert proved that a non-negative real quartic form $f(x,y,z)$ is
the sum of three squares of quadratic forms. We give a new proof
which shows that if the plane curve $Q$ defined by $f$ is
smooth, then $f$ has exactly $8$ such representations, up to
equivalence. They correspond to those real $2$-torsion points of the
Jacobian of $Q$ which are not represented by a
conjugation-invariant divisor on $Q$.

\vskip 0.5\baselineskip

%\selectlanguage{francais}
%
%\noindent{\bf R\'esum\'e}
%
%\noindent
%Hilbert a d\'emontr\'e qu'une forme r\'eelle non n\'egative $f(x,y,z)$
%de degr\'e 4 est la somme de trois carr\'es de formes quadratiques. Nous
%donnons une nouvelle d\'emonstration qui montre que si la courbe plane
%$Q$ definie par $f$ est non singuli\`ere, alors $f$ a exactement $8$
%telles repr\'esentations, \`a equivalence pr\`es. Elles correspondent
%aux points de $2$- torsion du jacobien de $Q$ qui ne sont pas
%repr\'esent\'es par un diviseur de $Q$ invariant par conjugaison.
%
%
%\vskip 0.5\baselineskip
%\noindent
\end{abstract}

\end{frontmatter}
%\mbox{\ }\vspace{-45pt}

%-------------------------------------------------------------------%

\selectlanguage{english}
\section{Introduction}

A \emph{ternary quartic} is a homogeneous polynomial $f(x,y,z)$
of degree~$4$ in three variables. If $f$ has real coefficients, then
$f$ is \emph{non-negative} if $f(x,y,z)\ge0$ for all real $x$,
$y$, $z$. Hilbert~\cite{Hi1888} showed that every non-negative real
ternary quartic form is a sum of three squares of quadratic forms.
His proof (see \cite{Rudin}, \cite{Swan} for modern expositions) was
non-constructive: The map
$$\pi\colon\ (p,q,r)\ \longmapsto\ p^2+q^2+r^2$$
from triples of real quadratic forms to non-negative quartic forms is
surjective, as it is both open and closed when restricted to the
preimage of the (dense) connected set of non-negative quartic forms
which define a smooth complex plane curve.  An elementary and
constructive approach to Hilbert's theorem was recently begun by
Pfister \cite{Pf04}.

A \emph{quadratic representation} of a complex ternary quartic form
$f=f(x,y,z)$ is an expression
 \begin{equation}\label{E:Q-rep}
   f\ =\ p^2 + q^2 + r^2\
 \end{equation}
where $p$, $q$, $r$ are complex quadratic forms. 
A representation $f=(p')^2 + (q')^2 + (r')^2$ is \emph{equivalent} 
to this if $p,q,r$ and $p',q',r'$ have the same linear
span in the space of quadratic forms.

Powers and Reznick \cite{PR00} investigated quadratic representations
computationally, using the Gram matrix method of~\cite{CLR}. In
several examples of non-negative ternary quartics, they always found
$63$ inequivalent representations as a sum of three squares of
complex quadratic forms; $15$ of these were sums or differences of
squares of real forms. We explain these numbers, in particular the
number $15$, and show that precisely $8$ of the $15$ are sums of
squares.

If the complex plane curve $Q$ defined by $f=0$ is
smooth, it has genus $3$, and so the Jacobian $J$ of
$Q$ has $2^6-1=63$ non-zero $2$-torsion points. Coble~\cite{Co1929}
showed that these are in one-to-one correspondence with equivalence
classes of quadratic representations of $f$.
If $f$ is real, then $Q$ and $J$ are defined over $\R$. The non-zero
$2$-torsion points of $J(\R)$ correspond to \emph{signed
quadratic representations}
$f=\pm p_1\pm p_2\pm p_3$, where 
$p_i$ are real
quadratic forms. If $f$ is also non-negative, the real Lie group
$J(\R)$ has two connected components, and hence has $2^4-1=15$
non-zero $2$-torsion points. 
We use Galois cohomology to determine which $2$-torsion points give rise to sum
of squares representations over $\R$.

\begin{theorem}\label{T:main}
Suppose that $f(x,y,z)$ is a non-negative real quartic form which defines a
smooth complex plane curve $Q$. Then the inequivalent
representations of $f$ as a sum of three squares are in one-to-one
correspondence with the eight $2$-torsion points in the non-identity
component of $J(\R)$, where $J$ is the Jacobian of $Q$.
\end{theorem}

Wall~\cite{Wa91} studies quadratic representations of 
(possibly singular) complex ternary quartic forms $f$. Again, in the 
irreducible case, the
non-trivial $2$-torsion points on the generalized Jacobian give equivalence
classes of quadratic representations of $f$. These representations
are special in that they have no basepoints.

Quadratic representations with a given base locus $B$ 
correspond to the $2$-torsion points on the Jacobian of a curve  
$\wt Q$, which is the image of $Q$ under the complete linear series
of quadrics through $B$. Classifying all possibilities for $B$ 
gives the number of inequivalent quadratic representations of
$f$. If $f$ is real and
non-negative, this classification, together with arguments from
Galois cohomology, gives all inequivalent representations of $f$ as a
sum of squares. 
This complete analysis will appear in an
unabridged version.

We thank C.T.C.~Wall, who brought his work to our attention, and the
organizers of the RAAG conference in Rennes in June 2001, where this
work began.

%-------------------------------------------------------------------%

\section{Basepoint-free quadratic representations}\label{S:PFQR}

Let $f(x,y,z)$ be an irreducible quartic form over $\C$, and let $Q$
be the complex plane curve $f=0$. The Picard group
$\Pic(Q)$ of $Q$ is the group of Weil divisors on the regular part of
$Q$, modulo divisors of rational functions which are
invertible around the singular locus of $Q$.
Let $J_Q$ be the generalized Jacobian of $Q$, so that $J_Q(\C)$ is
the identity component of $\Pic(Q)$.
Its structure may be determined from the Jacobian of the normalization
$\wt Q$ of $Q$ via the exact sequence \cite[Ex.~II.6.9]{Ha77}
 $$
   0\ \longrightarrow\
    \bigoplus_{p\in Q} \wt\calO^*_p/\calO^*_p
    \ \longrightarrow\ J_Q(\C)\ \longrightarrow\
    J_{\wt Q}(\C)\ \longrightarrow\ 0\,,
 $$
where $\calO_p$ is the local ring of $Q$ at $p$, $\wt\calO_p$ is its
normalization, and $^*$ indicates the group of units.

The \emph{base locus} $B$  of a quadratic representation~\eqref{E:Q-rep}
of $f$ is the
zero scheme of the homogeneous ideal generated by the forms $p$, $q$,
$r$. The closed subscheme $B$ is supported on the singular locus of
$Q$. We say that \eqref{E:Q-rep} is \emph{basepoint-free} if $B$ is
empty.

\begin{prop}[Coble~\cite{Co1929}, Wall~\cite{Wa91}]
The non-trivial $2$-torsion points of $J_Q$ are in one-to-one
correspondence with the equivalence classes of basepoint-free
quadratic representations of $f$.
\end{prop}

\noindent{\it Proof.\/}
Given a basepoint-free quadratic representation~\eqref{E:Q-rep},
consider the map
$$\varphi\colon\P^2\to\P^2,\quad x\mapsto\bigl(p(x):q(x):r(x)
\bigr).$$
The image of $Q$ under $\varphi$ is the conic $C$ defined by the
equation $y_0^2+y_1^2+y_2^2=0$. Let $y$ be a point in $C$ whose
preimages are regular points of $Q$. Then $\varphi^*(y)$ is an
effective divisor of degree $4$ that is not the divisor of a linear
form. Indeed, after a linear change of coordinates we can assume $y=
(0:1:i)$. 
A linear form vanishing on $\varphi^*(y)$
would divide each conic $\alpha p+\beta(q+ir)$ through $\varphi^*(y)$, and thus
would divide 
 $$
  f\ =\ p^2+(q+ir)(q-ir)\,,
 $$
contradicting the irreducibility of $f$.

Fix a linear form $\ell$ that does not vanish at any singular point
of $Q$. Then $L:=\div(\ell)$ is an effective divisor of degree~$4$ on
$Q$. Let $\zeta=[\varphi^*(y)-L]$. Since $2y$ is the divisor of a
linear form (the tangent line to $C$ at $y$), $\varphi^*(2y)$ is the
divisor on $Q$ of a quadratic form. Thus $2\zeta=0$. Moreover, $\zeta
\ne0$ as $\varphi^*(y)$ is not the divisor of a linear form. The
$2$-torsion point $\zeta$ of $J_Q$ depends only upon the map $\varphi$.

Conversely, suppose that $\zeta\in J_Q(\C)$ is a non-zero $2$-torsion
point. Let $D\ne D'$ be effective divisors which represent the class
$\zeta+[L]$ in $\Pic(Q)$. As $Q$ has arithmetic genus 3, the
Riemann-Roch Theorem implies that there is a pencil of such divisors.
%
%   Hartshorne, IV.1.9(a) explains this:
%
%   Riemann-Roch for singular curves implies that
%
%  dim H0(D) >= deg D + 1 - p_a  (p_a := arithmetic degree)
%
% Here,
%
%  deg D = 4 and p_a = 3, so that  dim H0(D) >= 2
%
Then $2D$, $2D'$ and $D+D'$ are effective divisors of degree~$8$, and
are all linearly equivalent to $2L$, the divisor of a conic. By the
Riemann-Roch Theorem there are
%
% Same source.  Let K be the divisor of the dualizing sheaf of Q.  This has
% degree 8, so that h0(K-2L)=0 (the divisor has negative degree).
%
%  Thus
%
%    h0(2L)-h0(K-2L) = deg 2L + 1 - p_a = 6
%
% Since there is a 6-dimensional space of quadrics, these cut out all divisors
% equivalent to 2L, as long as no two degree 8 divisors lie on the same quadric.
% Bezout's Theorem ensures this is the case.
%
quadratic forms  $q_0$, $q_1$ and $q_2$ such that
\begin{center}
$\div(q_0)=2D$, \quad $\div(q_1)=2D'$ \quad and \quad $\div(q_2)=D+D'$.
\end{center}

Therefore, the rational function $g:=q_0q_1/q_2^2$ on $Q$ is
constant. 
Scaling $q_1$ and $q_2$ appropriately, we may assume that 
$g\equiv1$ on $Q$ and also that $f=q_0q_1-q_2^2$.
Diagonalizing the quadratic form $q_0q_1-q_2^2$ 
gives a quadratic representation
for $f$. This defines the inverse of the previous map.
\qed

%-------------------------------------------------------------------%

\section{Quadratic representations of real quartics}

Suppose now that $f$ is a non-negative real quartic form defining a real plane
curve $Q$ with complexification $Q_\C=Q\otimes_\R\C$. 
The elements of $\Pic (Q)$ can be identified with those
divisor classes in $\Pic(Q_\C)$ that are represented by a conjugation-invariant
divisor. Let $J$ be the generalized Jacobian of $Q$.

If $\zeta\in J(\C)$ is the $2$-torsion point corresponding to a
signed quadratic representation
 $$
   f\ =\ \pm p^2\pm q^2 \pm r^2
 $$
consisting of real polynomials $p$, $q$, $r$, then $\zeta=\ol\zeta$,
i.e., $\zeta\in J(\R)$.

Conversely, let $0\ne\zeta\in J(\R)$ with $2\zeta=0$. Choose a real
linear form $\ell$ not vanishing on the singular points of $Q$, and
let $L=\div(\ell)$. %, an effective divisor on $Q$ of degree~$4$.
We can choose effective divisors $D\ne\ol D$ on
$Q_\C$ representing the class $\zeta+[L]$.
Then $2D$, $2\ol D$ and $D+\ol D$ are each
equivalent to $2L$. Let $r$ be a real quadratic form
%
%  By Riemann-Roch, any such quadratic form is unique, up to a scalar.
%
with divisor $D+\ol D$, and let $g$ be a (complex) quadratic form
with divisor $2D$ (both divisors taken on $Q_\C$). 
%Write $g=p+iq$, where $p$ and $q$ are the real and imaginary parts.

Since $D\sim\ol D$, there is a rational function $h$ on $Q_\C$,
invertible around $Q_\sing$, with  $\div(h)=D-\ol D$. Let $c=h\ol h$,
a nonzero real constant on $Q$. Since $\div(r)=\div(g)+\div(h)$,
there is a complex number $\alpha\ne0$ with $\frac rg=\alpha h$ on $Q$,
which implies that
 $$
   c|\alpha|^2\ =\ \frac rg\>\frac{\ol r}{\ol g}\ =\ \frac{r^2}{p^2+q^2}
 $$
on $Q$, where $p,q$ are the real and imaginary parts of 
$g=p+iq$. So the quartic form
$$
  u\ :=\ r^2-c\,|\alpha|^2(p^2+q^2)
$$
vanishes identically on $Q$. Since $u\ne0$,
% as we saw before in the complex case
$f$ is a constant
multiple of $u$. If $c>0$, we get a signed quadratic representation of $f$,
with both signs $\pm$ occuring.
If $c<0$, $f$ must be a positive multiple of $u$ since $f$
is non-negative, and we get a representation of $f$ as a sum of three
squares of real forms.

We now calculate the sign of $c$. For this we use the exact sequence
 \begin{equation}\label{HS}
   0\to\Pic(Q)\to\Pic(Q_\C)^G\labelto\partial\Br(\R)\to H^2_\et
  (Q,\G_m)
 \end{equation}
of \'etale cohomology groups. It arises from the Hochschild-Serre
spectral sequence for the Galois covering $Q_\C\to Q$ and coefficients
$\G_m$. Here $G=\Gal(\C/\R)$ acts on $\Pic(Q_\C)$ by conjugation, and
$\Pic(Q_\C)^G$ is the group of $G$-invariant divisor classes.
Moreover, $\Br(\R)=H^2_\et(\Spec\R,\G_m)$ is the Brauer group of
$\R$ (which is of order~$2$), and $\Br(\R)\to H^2_\et(Q,\G_m)$ is
the restriction map.

It is easy to see that
%%%%%%%%%%%%%%%%%%%%%%%%%%
% $\partial(\zeta)=(-1,c)$.
%
%   Frank does not see that $\partial(\zeta)=(-1,c)$ is so helpful.
%  (Is it that  H^2_\et(\Spec\R,\G_m) is the space of G_m-line bundles 
%   on \R ?)   Anyway, this is a little shorter and no more obfuscatory.
%
%In other words,
%%%%%%%%%%%%%%%%%%%%%%%%%%%%%%%%%%%%%%%%%%%5
$c<0$ if and only if $\partial(\zeta)$ is the non-trivial class in
$\Br(\R)$.
If $Q$ has an $\R$-point, then $\Br(\R)\to H^2_\et(Q,\G_m)$ has a
splitting given by that point, and hence $\partial$ vanishes
identically.

If $Q$ is smooth, then $f$ non-negative forces $Q(\R)=\emptyset$, 
and the map $\Br(\R)\to H^2_\et(Q,\G_m)$ is zero. 
In this case, $Pic(Q_C)^G$ contains an odd degree divisor if and only 
if the genus of $Q$ is even and $J(\R)^0$, the identity connected
component of the real Lie group $J(\R)$, is the kernel of the
restriction $J(\R)\to\Br(\R)$ of $\partial$~\cite{We1883,Ge66}. 
%In this case it is classically known (Weichold~\cite{We1883}, 
%Geyer~\cite{Ge66}) that $Pic(Q_C)^G$ contains an odd degree divisor if and only 
%if the genus of $Q$ is even, and that $J(\R)^0$, the identity connected
%component of the real Lie group $J(\R)$, is the kernel of the
%restriction $J(\R)\to\Br(\R)$ of $\partial$. 
Since in our case $g(Q)=3$, this implies that the sequence
$$0\to J(\R)^0\to J(\R)\labelto\partial\Br(\R)\to0$$
is (split) exact. If $Q$ is singular with
$Q(\R)=\emptyset$, one compares sequence \eqref{HS} for $Q$ to the
same sequence for the normalization $\wt Q$ of $Q$ and concludes that
$\partial\colon J(\R)\to\Br(\R)$ is surjective as well.
% (!!)

We complete the proof of Theorem \ref{T:main}. Since $f$ is
non-negative and $Q$ is smooth of genus $3$, 
$J(\R)^0\cong (S^1)^3$ as a real Lie group. By the facts just
mentioned, there exist $2^4-1=15$ non-zero $2$-torsion elements in
$J(\R)$.  The $8$ that do not lie in $J(\R)^0$, or
equivalently, which cannot be represented by a conjugation-invariant
divisor on $Q_\C$, are precisely those that give
rise to sums of squares representations of $f$.

%-------------------------------------------------------------------%

\providecommand{\bysame}{\leavevmode\hbox to3em{\hrulefill}\thinspace}
\providecommand{\MR}{\relax\ifhmode\unskip\space\fi MR }
\providecommand{\MRhref}[2]{%
  \href{http://www.ams.org/mathscinet-getitem?mr=#1}{#2}
}
\providecommand{\href}[2]{#2}

\end{document}